\documentclass[letterpaper,12pt]{amsart}
\usepackage[margin=1.5in]{geometry}
\usepackage{geometry}
\usepackage{amssymb}
\usepackage{amsmath}
\usepackage{xypic}
\usepackage{graphicx}
\usepackage{mathrsfs}
\usepackage{enumerate}
\usepackage[mathscr]{eucal}
\usepackage{setspace}
\usepackage{blkarray}
\usepackage{multicol}
\usepackage{lscape}
\usepackage{mathtools}
\usepackage[noadjust]{cite}

\theoremstyle{plain}
\newtheorem{theorem}{Theorem}
\newtheorem{lem}{Lemma}

\newtheorem{cor}{Corollary}
\newtheorem*{claim}{Claim}

\theoremstyle{definition}
\newtheorem*{definition}{Definition}
\newtheorem{conjecture}{Conjecture}

\theoremstyle{remark}
\newtheorem*{rem}{Remark}

\newtheorem{notation}{Notation}

\newcommand{\beas}{\begin{eqnarray*}}
\newcommand{\eeas}{\end{eqnarray*}}
\newcommand{\bes} {\begin{equation*}}
\newcommand{\ees} {\end{equation*}}
\newcommand{\be} {\begin{equation}}
\newcommand{\ee} {\end{equation}}
\newcommand{\bea} {\begin{eqnarray}}
\newcommand{\eea} {\end{eqnarray}}

\DeclareMathOperator{\Spec}{Spec}

\DeclareMathOperator{\Rad}{Rad}

\DeclareMathOperator{\mm}{\mathfrak{m}}
\DeclareMathOperator{\A}{\mathfrak{A}}

\newcommand{\Aa}{\alpha}
\newcommand{\Bb}{\beta}
\newcommand{\Gg}{\gamma}
\newcommand{\Dd}{\delta}
\newcommand{\La}{\lambda}

\begin{document}

\bibliographystyle{amsalpha}       
\setlength{\itemsep}{14pt}

\title{Nearly commuting matrices}
\author{Zhibek Kadyrsizova}
\maketitle
\begin{abstract}
We prove that the algebraic set of pairs of matrices with a diagonal commutator over a field of positive prime characteristic, its irreducible components, and their intersection  are $F$-pure when the size of matrices is equal to 3. Furthermore, we show that this algebraic set is reduced and the intersection of its irreducible components is irreducible in any characteristic for pairs of matrices of any size. In addition, we discuss various conjectures on the singularities of these algebraic sets and the system of parameters on the corresponding coordinate rings. 
\end{abstract}

\section{Introduction and preliminaries}
In this paper we study algebraic sets of pairs of matrices such that  their commutator is either nonzero diagonal or zero. We also consider some other related algebraic sets.  First let us define relevant notions.

Let $X=\left(x_{ij}\right)_{1\leq i, j\leq n}$ and $Y=(y_{ij})_{1\leq i, j\leq n}$ be $n\times n$ matrices of indeterminates over a field $K$. Let $R=K[X,Y]$ be the polynomial ring in $\{x_{ij}, y_{ij}\}_{1\leq i, j\leq n}$ and let $I$ denote the ideal generated by the off-diagonal entries of the commutator matrix $XY-YX$ and  $J$ denote the ideal generated by the entries of $XY-YX$. The ideal $I$ defines the algebraic set of pairs of  matrices with a diagonal commutator and is called \emph{the algebraic set of nearly commuting matrices}. The ideal $J$ defines the algebraic set of pairs of commuting matrices.  

Let $u_{ij}$ denote the $(i,j)$th entry of the matrix $XY-YX$. Then $I=(u_{ij} \, | \, 1\leq i\neq j \leq n)$ and $J=(u_{ij} \, | \,1\leq i, \, j \leq n)$. 

\begin{theorem}[\cite{G61}] The algebraic set of commuting matrices is irreducible, i.e., it is a variety. Equivalently, $\Rad(J)$ is prime. \end{theorem}

The following results are due to A. Knutson \cite{KNU05}, when the characteristic of the field is 0, and to H.Young \cite{Y11} in all characteristics. 

\begin{theorem}  [\cite{KNU05}, \cite{Y11}] The algebraic set of nearly commuting matrices is a complete intersection, with the variety of commuting matrices as one of its irreducible components.  In particular, the set $\{u_{ij}| 1\leq i\neq j\leq n\}$ is a regular sequence and the dimension of $R/I$ is $n^2+n$.  \end{theorem}

\begin{theorem} [\cite{KNU05}, \cite{Y11}] When $K$ has characteristic zero, $I$ is a radical ideal. \end{theorem}

A. Knutson in his paper \cite{KNU05} conjectured that $\mathbb{V}(I)$ has only two irreducible components and it was proved in all characteristics by H.Young in his thesis, ~\cite{Y11}.  

\begin{theorem} [\cite{Y11}] If $n \geq 2$, the algebraic set of nearly commuting matrices has two irreducible components, one of which is the variety of  commuting matrices and the other is the so-called skew component. That is, $I$ has two minimal primes, one of which is $\Rad(J)$.\end{theorem}

Let  $P=\Rad(J)$ and let  $Q$ denote the other minimal prime of $I$, i.e., $\Rad(I)=P\bigcap Q$.

The following conjecture was made in 1982 by M. Artin and  M. Hochster. 

\begin{conjecture} $J$ is reduced, i.e.,  $J=P$. \end{conjecture}

It was answered positively by Mary Thompson in her thesis in the case of $3\times 3$ matrices.  
\begin{theorem}[\cite{T85}] $R/J$ is a Cohen-Macaulay domain when $n=3$. \end{theorem}
Now let us go back to algebraic sets of nearly commuting matrices and their irreducible components. First, we take a look at what we have when $n=1,2$. 

When $n=1$, everything is trivial. More precisely, $I=P=Q=K[x_{11}, y_{11}]$.

When $n=2$, without loss of generality we may replace $X$ and $Y$ by $X-x_{22}\text{I}_n$ and $Y-y_{22}\text{I}_n$ respectively. Here $\text{I}_n$ is the identity matrix of size $n$.  Denote $x_{11}'=x_{11}-x_{22}, \, y_{11}'=y_{11}-y_{22}$. Then the generators of $I$ are 2 by 2 minors

$$u_{12}=\left | \begin{array}{cc}
                           x_{11}' & x_{12} \\
                           y_{11}' &y_{12}  \end{array} \right|,  \quad
u_{21}=-\left | \begin{array}{cc}
                           x_{11}' & x_{21} \\
                           y_{11}' &y_{21}  \end{array} \right|. $$
The diagonal entries of $XY-YX$ are      
$$u_{11}=\left | \begin{array}{cc}
                           x_{12} & x_{21} \\
                           y_{12} &y_{21}  \end{array} \right|,  \quad
u_{22}=-\left | \begin{array}{cc}
                           x_{12} & x_{21} \\
                           y_{12} &y_{21}  \end{array} \right|.$$

Then  $J$ is the ideal generated by size 2 minors of $\left [ \begin{array}{ccc}
                           x_{11}' & x_{12} & x_{21} \\
                           y_{11}' &y_{12} &y_{21}  \end{array} \right]$ and therefore, $J=P$ is prime.  We also have that $Q=(x_{11}', y_{11}')$. Moreover, $I=P\bigcap Q$ is radical and $P+Q=(x_{12}y_{21}-x_{21}y_{12}, x'_{11}, y'_{11})$ is prime. 
                           
We have that $$(u_{12}u_{21})^{p-1}=(x_{11}'y_{12}-x_{12}y_{11}')^{p-1}(x_{11}'y_{21}-x_{21}y_{11}')^{p-1}=$$
$$\sum_{\Aa=0}^{p-1}\sum_{\Bb=0}^{p-1}(-1)^{\Aa+\Bb}{p-1 \choose \Aa}{p-1 \choose \Bb} (x_{11}')^{\Aa+\Bb}(y_{11}')^{2(p-1)-\Aa-\Bb}x_{12}^{p-1-\Aa}y_{12}^{\Aa}x_{21}^{p-1-\Bb}y_{21}^{\Bb}.$$
Therefore,  $(u_{12}u_{21})^{p-1}$ has a monomial term $(x_{11}'y_{11}'x_{12}y_{21})^{p-1}$ with coefficient $(-1)^{p-1}$.   Since $I^{[p]}:I=(u_{12}u_{21})^{p-1}+I^{[p]}$, $R/I$ is $F$-pure, see Fedder's criterion Lemma ~\ref{Fed2}.   Furthermore, determinantal rings $R/P$, $R/Q$, $R/(P+Q)$ are $F$-regular, see \cite{HH94}. 

Therefore, for the rest of the paper we shall use the following notations.
\begin{notation} \label{NotComm}
Let $n \geq 3$ be an integer. Let $X=(x_{ij})_{1\leq i, j\leq n}$ and $Y=(y_{ij})_{1\leq i, j\leq n}$ be $n\times n$ matrices of indeterminates over a field $K$. Let $R=K[X,Y]$ be the polynomial ring in $\{x_{ij}, y_{ij}\}_{1\leq i, j\leq n}$ and let $I$ denote the ideal generated by the off-diagonal entries of the commutator matrix $XY-YX$ and  $J$ denote the ideal generated by the entries of $XY-YX$. Let $P$ denote the radical of $J$ and $Q$ be the other minimal prime of $\Rad(I)$.
\end{notation}
We prove the following results in this paper.
\begin{theorem} Let $R$ be a ring as in Notation ~\ref{NotComm}.  Assume also that the field $K$ has  positive prime characteristic. Then $R/I$, $R/P$, $R/Q$ and $R/(P+Q)$ are $F$-pure rings when $n=3$. In other words, the algebraic set of nearly commuting matrices of size 3, its irreducible components and their intersection are $F$-pure. In particular, the skew component is reduced in this case. \end{theorem}
\begin{theorem} Let $R$ be a ring as in Notation ~\ref{NotComm}.  Then $R/I$ is reduced. In other words, the algebraic set of nearly commuting matrices is reduced for matrices of all sizes and in all characteristics.  \end{theorem}
\begin{theorem}  The intersection of the variety of commuting matrices and the skew component is irreducible, that is, $\Rad(P+Q)$ is prime.  \end{theorem}

\section{$F$-purity}
In this section we show that the coordinate ring of the algebraic set of pairs of matrices with a diagonal commutator is $F$-pure in the case of 3 by 3 matrices. Moreover, we also show that it implies the corresponding fact for its irreducible components, the variety of commuting matrices and the skew-component,  and their intersection. 

First we state two lemmas due to R. Fedder and they include a criterion for $F$-purity for finitely generated $K$-algebras and which has a particularly convenient form  for complete  intersections. 
\begin{lem}[Fedder ~\cite{F87}] \label{Fed1}
Let $S$ be a regular local  ring or a polynomial ring over a field. If $S$ has characteristic $p>0$ and $I$ is an unmixed proper ideal (homogeneous in the polynomial case) with the primary decomposition $I=\bigcap_{i=1}^n \A_i$, then $I^{[p]}:I=\bigcap_{i=1}^n (\A^{[p]}:\A)$. 
\end{lem}
\begin{lem}[Fedder's criterion ~\cite{F87}] \label{Fed2} Let $(S, \mm)$ be a regular local  ring or a polynomial ring over a field with its (homogeneous) maximal ideal. If $S$ has characteristic $p>0$ and  $I$ is a proper ideal (homogeneous in the polynomial case),  then $S/I$ is $F$-pure if and only if $I^{[p]}:I \not \subset \mm^{[p]}$. \end{lem}
The next result is a straightforward consequence of the above two lemmas. It will prove to be quite useful for us. 
\begin{lem} \label{Cool} Let $S$ be a regular local ring or a polynomial ring over a field. Suppose that $S$ has characteristic $p>0$ and $I$ is an ideal of $S$ (homogeneous in the polynomial case). Suppose also that $S/I$ is $F$-pure and $I=\bigcap_{i=1}^n \A_i$ is the primary decomposition. Then $S/(\A_{i_1}+\ldots \A_{i_m})$ is  $F$-pure for all $1\leq i_1< ...<i_m\leq n$ and for all $1\leq m \leq n$. 
\end{lem}
\proof  Observe first that $(\A_{i_1}+\ldots \A_{i_m})^{[p]}:_S (\A_{i_1}+\ldots \A_{i_m})\supseteq \bigcap_{j=1}^m ( \A_{i_j}^{[p]}: \A_{i_j})\supseteq  \bigcap_{i=1}^n ( \A_{i}^{[p]}: \A_{i})=(I^{[p]}:I)$. The rest is immediate from Lemma ~\ref{Fed1} and Lemma ~\ref{Fed2}. \qed

The above lemma is closely related to results on compatibly split ideals, cf. \cite{ST12}. 

Immediately we get the corresponding result for our algebraic set. 

\begin{cor} Suppose that the coordinate ring of the algebraic set of  nearly commuting matrices $R/I$ is $F$-pure. Then $R/P$, $R/Q$ and $R/(P+Q)$ are $F$-pure.
\end{cor}

Next, we use Fedder's criterion to show $F$-purity of $R/I$ in case when $n=3$. 
\begin{theorem} \label{ThPur}Let $K$ be a field of characteristic $p>0$ and let $n=3$. Let $R$ be a ring as in Notation ~\ref{NotComm}. Then $R/I$ is $F$-pure. 
\end{theorem}

\begin{proof}  Recall that $I$ is generated by a regular sequence $\{u_{ij}| 1\leq i\neq j \leq n\}$. Therefore, $I^{[p]}:I=\left(\prod_{1\leq i\neq j \leq n}u_{ij}^{p-1}\right)R+I^{[p]}$. Thus by Fedder's criterion it is sufficient to prove that $\prod_{1\leq i\neq j \leq n}u_{ij}^{p-1} \notin \mm^{[p]}$, where $\mm$ is the homogeneous maximal ideal in $R$. We show this by proving the following claim. 

\begin{claim} If $\mu=x_{12}x_{13}x_{21}x_{23}x_{31}x_{33}y_{11}y_{12}y_{23}y_{31}y_{32}y_{33}$, then $\mu^{p-1}$ is a monomial term of $\prod_{1\leq i \neq j \leq 3}u_{ij}^{p-1}$ with a nonzero coefficient modulo $p$. 
\end{claim}
\proof 

We compute the coefficient of $\mu^{p-1}$. It can be obtained by choosing a monomial from every $u_{ij}$ in the following way:\\
$u_{12}: \quad (-x_{12}y_{11})^{\alpha_1}(x_{13}y_{32})^{\beta_1}$\\
$u_{13}: \quad (-x_{23}y_{12})^{\alpha_2}(x_{12}y_{23})^{\beta_2}(-x_{13}y_{11})^{\gamma_2}(x_{13}y_{33})^{\delta_2}$\\
$u_{21}: \quad (-x_{31}y_{23})^{\alpha_3}(x_{21}y_{11})^{\beta_3}(x_{23}y_{31})^{\gamma_3}$\\
$u_{23}: \quad (x_{23}y_{33})^{\alpha_4}(-x_{33}y_{23})^{\beta_4}$\\
$u_{31}: \quad (-x_{21}y_{32})^{\alpha_5}(x_{33}y_{31})^{\beta_5}(-x_{31}y_{33})^{\gamma_5}(x_{31}y_{11})^{\delta_5}$\\
$u_{32}: \quad (x_{31}y_{12})^{\alpha_6}(-x_{12}y_{31})^{\beta_6}(x_{33}y_{32})^{\gamma_6}$\\
Then the exponents $A_{st}$ and $B_{st} $ of each $x_{st}$ and $y_{st}$ respectively are 
\begin{multicols}{2}
$A_{12}=\alpha_1+\beta_2+\beta_6$

$A_{13 }=\beta_1+\gamma_2+\delta_2$

$A_{21 }=\beta_3+\Aa_5$

$A_{23 }=\alpha_2+\gamma_3+\alpha_4$

$A_{31}=\alpha_3+\gamma_5+\delta_5+\alpha_6$

$A_{33}=\beta_4+\beta_5+\gamma_6$
\end{multicols}

\begin{multicols}{2}
$B_{11}=\alpha_1+\gamma_2+\beta_3+\Dd_5$

$B_{12 }=\alpha_2+\alpha_6$

$B_{23 }=\beta_2+\Aa_3+\Bb_4$

$B_{31 }=\gamma_3+\Bb_5+\Bb_6$

$B_{ 32}=\Bb_1+\Aa_5+\Gg_6$

$B_{33 }=\delta_2+\Aa_4+\Gg_5$
\end{multicols}

In addition, denote
\begin{multicols}{2}
 $C_{12}=\Aa_1 +\Bb_1$, 
 
 $C_{13}=\Aa_2 +\Bb_2 +\Gg_2 + \Dd_2$,
 
 $ C_{21}=\Aa_3 +\Bb_3 +\Gg_3$, 
 
 $C_{23}=\Aa_4 +\Bb_4$, 
 
 $C_{31}=\Aa_5 +\Bb_5 +\Gg _5+ \Dd_5$, 
 
 $C_{32}=\Aa_6 +\Bb_6 +\Gg_6$.
  \end{multicols}
Our goal is to find all nonnegative integer tuples $\Aa=(\Aa_1, \ldots, \Aa_6), \, \Bb=(\Bb_1, \ldots, \Bb_6), \Gg=(\Gg_2,\Gg_3,\Gg_5, \Gg_6), \, \Dd=(\Dd_5, \Dd_6) $ such that $A_{st}=p-1$, $B_{st}=p-1$ for all $1\leq s,t \leq 3 $ and $C_{ij}=p-1$ for all $1\leq i\neq j\leq 3$. 

Notice that the linear system does not have a nonzero determinant: the sum of the first 12 equations is twice the sum of the rest 6 equations. Therefore,  there is not a unique solution. 

The above linear system can be solved using standard methods from linear algebra and has the following solution:

\begin{small}
$$ \left[ \begin{array}{c}
                    \alpha \\
                    \beta\\
                    \gamma\\
                    \delta\\
                     \end{array} \right]=$$$$\left[ \begin{array}{cccccc}
                    a        & b    &d     & p-1-b-a+d       &a        &p-1-b\\
                    p-1-a & p-1-a-b & p-1-a & a+b-d& p-1-a-b+d &b\\
                   -   &a-c&a-d&-     &b+a-d-c  &0\\
                   - &  c  &- & -   &  c-a  & - \\
                        \end{array} \right]$$
\end{small}                        
where the column vector $[\Aa, \Bb, \Gg, \Dd]$ represents the  matrix of solutions and $a, b, c,d$ are arbitrary non-negative integers. 

Since we look for non-negative integer solutions we must have that $a=c$ and $a,\,  b \geq d$ and $a+b\leq p-1$.  Hence we have that

\begin{small}

$$ \left[ \begin{array}{c}
                    \alpha \\
                    \beta\\
                    \gamma\\
                    \delta\\
                     \end{array} \right]=$$$$\left[ \begin{array}{cccccc}
                    a        & b    &d     & p-1-b-a+d       &a        &p-1-b\\
                    p-1-a & p-1-a-b & p-1-a & a+b-d& p-1-a-b+d &b\\
                   -   &0&a-d&-     &b-d  &0\\
                   - &  a  &- & -   &  0  & - \\
                        \end{array} \right]$$
\end{small}

Therefore, the coefficient of $\mu^{p-1}$ is the sum of expressions of the form 
$$(-1)^{\Aa_1+\Aa_2+\Gg_2+\Aa_3+\Bb_4+\Aa_5+\Gg_5+\Bb_6}((p-1)!)^6/(\Aa_1! \ldots \Aa_6! \Bb_1! \ldots \Bb_6!\Gg_2!\Gg_3!\Gg_5! \Gg_6!\Dd_5! \Dd_6!)$$
where $\Aa=(\Aa_1, \ldots, \Aa_6), \, \Bb=(\Bb_1, \ldots, \Bb_6), \Gg=(\Gg_2,\Gg_3,\Gg_5, \Gg_6), \, \Dd=(\Dd_5, \Dd_6)$ run over all solutions of the linear system above. That is, 
$$\sum_{d=0}^{{(p-1)/2}}\sum_{a,b \geq d, \, a+b\leq p-1} (-1)^{a-d}{p-1\choose  a}^2{p-1\choose  b}^2{p-1\choose  a+b-d}^2{p-1 -b \choose  a}{a+b-d\choose  b}{b\choose  d}$$
which modulo $p$ is equivalent to 
$$\sum_{d=0}^{{(p-1)/2}}\sum_{a,b \geq d, \, a+b\leq p-1} (-1)^{a-d}{p-1 -b \choose  a}{a+b-d\choose  b}{b\choose  d}$$
It also can be written as
$$\sum_{d=0}^{{(p-1)/2}}\sum_{a,b \geq d, \, a+b\leq p-1} (-1)^{a-d}{p-1 -b \choose  a}{a+b-d\choose  a-d \quad  b-d \quad  d}$$
or
$$\sum_{b=0}^{p-1}\sum_{d=0}^{b}\sum_{a=d}^{p-1-b} (-1)^{a-d}{p-1 -b \choose  a}{a+b-d\choose  b}{b\choose  d}$$
The following lemma shows that the above expression is equal to 1 for all values of $p$. In fact,  for this purpose $p$ does not have to be prime. 

\begin{lem} Let $C_m=\sum_{b=0}^m \sum_{d=0}^{b}\sum_{a=d}^{m-b} (-1)^{a-d}{m -b \choose  a}{a+b-d\choose  b}{b\choose  d}$. Then $C_m=1$ for all  $m\geq 1$. \end{lem}
\proof 
We prove a stronger statement.
\begin{claim}  
Let $B_{m, b}=\sum_{d=0}^b\sum_{a=d}^{m-b}(-1)^{a-d}{m -b \choose  a}{a+b-d\choose  b}{b\choose  d}$.   Then 
for all $m \geq 1$
$$ B_{m,b} =\left \{ \begin{array}{cc}
                                              0& \mbox{if $0 \leq b\leq m-1$};    \\
                                              1&  \mbox{if $b=m$}.  \\
                                               \end{array} \right .$$
\end{claim}
\proof 
First observe that $B_{m,m}=\sum_{d=0}^m\sum_{a=d}^{0}(-1)^{a-d}{m -b \choose  a}{a+b-d\choose  b}{b\choose  d}=1$ and $B_{m,0}=\sum_{a=0}^{m}(-1)^{a}{m \choose  a}=0$. Hence we may assume that $0<b <m$. 

Let $A_{m,b,d}=\sum_{a=d}^{m-b}(-1)^{a}{m -b \choose  a}{a+b-d\choose  b}$, then $B_{m,b}=\sum_{d=0}^b (-1)^{d}{b\choose  d}A_{m,b,d}$. Consider the  difference
$$A_{m,b,d}-A_{m,b, d+1}=$$
$$\sum_{a=d}^{m-b}(-1)^{a}{m -b \choose  a}{a+b-d\choose  b}-\sum_{a=d+1}^{m-b}(-1)^{a}{m -b \choose  a}{a+b-d-1\choose  b}=$$
$$(-1)^d{m-b \choose d}+\sum_{a=d+1}^{m-b}(-1)^{a}{m -b \choose  a}({a+b-d\choose  b}-{a+b-d-1\choose  b})=$$
Using Pascal's identity, we get 
$$(-1)^d{m-b \choose d}+\sum_{a=d+1}^{m-b}(-1)^{a}{m -b \choose  a}{a+b-d-1\choose  b-1}=$$
$$\sum_{a=d}^{m-b}(-1)^{a}{m -b \choose  a}{a+b-d-1\choose  b-1}=$$
$$\sum_{a=d}^{m-1-(b-1)}(-1)^{a}{m-1 -(b-1) \choose  a}{a+(b-1)-d\choose  b-1}.$$
Thus we have that 
$$A_{m,b,d}-A_{m,b, d+1}=A_{m-1, b-1, d} \textnormal{   for all $m-1\geq b\geq d+1$ and $d\geq 0$}.$$
Therefore,
$$B_{m-1, b-1}=\sum_{d=0}^{b-1} (-1)^{d}{b-1\choose  d}A_{m-1,b-1,d}=\sum_{d=0}^{b-1} (-1)^{d}{b-1\choose  d}(A_{m,b,d}-A_{m,b, d+1})=$$
$$\sum_{d=0}^{b-1} (-1)^{d}{b-1\choose  d}A_{m,b,d}-\sum_{d=0}^{b-1} (-1)^{d}{b-1\choose  d}A_{m,b, d+1}=$$
$$\sum_{d=0}^{b-1} (-1)^{d}{b\choose  d}\frac{b-d}{b}A_{m,b,d}-\sum_{d=1}^{b} (-1)^{d-1}{b-1\choose  d-1}A_{m,b, d}=$$
$$\sum_{d=0}^{b-1} (-1)^{d}{b\choose  d}\frac{b-d}{b}A_{m,b,d}+\sum_{d=1}^{b} (-1)^{d}{b\choose  d}\frac{d}{b}A_{m,b, d}=$$
$$\sum_{d=1}^{b-1} (-1)^{d}{b\choose  d}A_{m,b,d}+A_{m,b, 0}+(-1)^bA_{m,b,b}=$$
$$\sum_{d=0}^{b} (-1)^{d}{b\choose  d}A_{m,b,d}=B_{m,b}.$$

Thus we have that $B_{m-1, b-1}=B_{m,b}$ for all $m \geq 1$ and $m-1\geq b\geq 1$.  

In case $m=1$, we only have  $B_{1,0}= 0$.
Finally, use induction on $m$ to conclude that $B_{m,b}=0$ for all $m\geq 1$ and $m-1\geq b$. 

Thus, $C_m=\sum_{b=0}^{m}B_{m,b}=1$. \qed

Finally,  we complete the proof of Theorem ~\ref{ThPur}.  We have that $\prod_{1\leq i\neq j \leq n}u_{ij}^{p-1} \notin \mm^{[p]}$ and $R/I$ is $F$-pure when $n=3$. 
\end{proof}
\begin{cor} Let $R$ be a ring as in Notation ~\ref{NotComm}.  When $n=3$, $R/P, R/Q$ and $R/(P+Q) $ are $F$-pure Cohen-Macaulay rings and $R/(P+Q)$ is Gorenstein. \end{cor}
\begin{proof}
By \cite{T85}, $R/P$ is a Cohen-Macaulay ring when $n=3$. Since  the ideals $P$ and $Q$ are linked via $I$, that is $I:P=Q$ and $I:Q=P$, we have that $R/Q$ is also Cohen-Macaulay, see \cite{PS74}. Moreover, the theory of linkage also implies that $(P+Q)/P$ and $(P+Q)/Q$ are isomorphic to the canonical modules of $R/P$ and $R/Q$, respectively. Hence $R/(P+Q)$ is Gorenstein of dimension $n^2+n-1$. 
\end{proof}
\begin{cor} Let $R$ be a ring as in Notation ~\ref{NotComm}. Then $P+Q$ is radical when $n=3$. \end{cor}

\begin{rem} We prove in the next section that  for all $n$, the radical of $P+Q$ is prime, which implies that $P+Q$ is prime when $n=3$. In particular, we have that $R/(P+Q)$ is a domain when $n=3$. \end{rem}
\vspace{5mm}
\section{Irreducibility of $P+Q$}
In this section we prove that the intersection of the variety of commuting matrices and the skew-component is irreducible. But first we define some notions.

\definition \label{poly} Let $X$ be an $n$ by $n$ matrix of indeterminates. Then $D(X)$ is an $n$ by $n$ matrix whose $i$th column is defined by the diagonal entries of $X^{i-1}$ numbered from upper left corner to lower right corner. Let $\mathcal{P}(X)$ denote the determinant of $D(X)$. 

\begin{theorem}[\cite{Y11}] $\mathcal{P}(X)$ is an irreducible polynomial. \end{theorem}

\begin{rem}  $\mathcal{P}(X)=\mathcal{P}(X-aI)$, where $a\in K$ and $I\in M_n(K)$ is the identity matrix. \end{rem}

The next two lemmas are due to H. Young. They give us the connection between the variety defined by $\mathcal{P}(X)$ and the algebraic set of nearly commuting matrices. 
\begin{lem}[\cite{Y11}] \label{L2Y} Given an $n \times n$ matrix $A$, if there exists a matrix $B$ such that $[A,B]$ is a non-zero diagonal matrix, then $\mathcal{P}(A)=0$. \end{lem}

\begin{lem}[\cite{Y11}] \label{L1Y} There is a dense open set  $\mathcal{U}$ in the variety defined by $\mathcal{P}(X)$ where for every point $A$ in $\mathcal{U}$, there exists a matrix $B$ such that $[A,B]$ is a nonzero diagonal matrix.  \qed \end{lem}

The following notion of a discriminant is of significant importance in matrix theory. We use it in this section in order to reduce our study to the case when commuting matrices have a particularly simple characterization. 
\begin{definition} Let $A \in M_n(K)$. Then the \emph{discriminant} $\Delta(A)$ of $A$ is the discriminant of its characteristic polynomial. That is, if $K$ contains all the eigenvalues $\La_1, \ldots, \La_n$ of $A$, then  $\Delta(A)=\prod_{1\leq i<j\leq n}(\La_i-\La_j)^2$.  \end{definition}

\textbf{Fact.} Let $A \in M_n(K)$ be a matrix such that $\Delta(A) \neq 0$, or equivalently, $A$ has distinct eigenvalues.  Then a matrix $B$ commutes with $A$ if and only if $B$ is a polynomial in A of degree at most $n-1$, see Theorem 3.2.4.2 \cite{Hor85}. 

\begin{rem} $\mathcal{P}(X)$ is an irreducible polynomial of degree $n(n-1)/2$ and $\Delta(X)$ is a polynomial of degree $n(n-1)$. Moreover, when $n \geq 3$,  $\mathcal{P}(X)$ does not divide $\Delta(X)$. This can be proved by showing that there exists a matrix $A$ with the property that $\mathcal{P}(A)=0$ while $\Delta(A) \neq 0$. For example, for this purpose one can use the following matrices.  
\begin{center}
$E_n=\left[ \begin{array}{cccccc}
0& 1&0& \ldots &0 &0\\
0& 0&1& \ldots &0 &0\\
& \ldots&&\ldots&&\\
0& 0&0& \ldots &0 &1\\
1& 0&0& \ldots &0 &0\\
 \end{array} \right] $  if $p\nmid n$, and  
 \vspace{4mm}
 $\widetilde{E}_n=\left[ \begin{array}{c|c}
                                                                              0 &0\\
                                                                              \hline
                                                                              0& E_{n-1}
                                                                           \end{array} \right]$, otherwise. 
\end{center}
The characteristic polynomials are $x^n-1$ for $E_n$ and $x(x^{n-1}-1)$ for $\widetilde{E}_n$. \end{rem}
Here is the outline for how we prove the main result of this section, that is,  $\Rad(P+Q)$ is prime.
\begin{enumerate}
\item[(1)] $\dim \mathbb{V}(P+Q)=n^2+n-1$ and  $\mathbb{V}(P+Q)$ is equidimensional.
\item[(2)] $ \dim \mathbb{V}(P+Q, \Delta(X))\leq n^2+n-2$, that is, $\Delta(X)$ is not in any of the minimal primes of $P+Q$.  
\item[(3)] $\Rad(P+Q)R_{\Delta(X)}$ is a prime ideal. 
\end{enumerate}

We observe that $P+Q$ has no minimal primes of height larger than one over $P$ and $Q$. First we need the following theorem due to R. Hartshorne. 
\begin{theorem}[\cite{Ha62} Proposition 2.1] Let $A$ be a Noetherian local ring with the maximal ideal $\mm$. If $\Spec(A) - \{\mm\}$ is disconnected, then the depth of $A$ is at most 1. 
\end{theorem}
\begin{lem} Let  $P$ and $Q$ be ideals as in Notation ~\ref{NotComm}. Then  every minimal prime of $P+Q$ has height $n^2-n+1$. \end{lem}
\begin{proof}
Suppose that there exists a minimal prime ideal $T$ of $P+Q$ of height at least $\text{ht}(I)+2$. 
Localize at $T$. Then  $(P+Q)(R/I)_T$ is $T(R/I)_T$-primary.  Moreover, $\mathbb{V}(P)$ and $\mathbb{V}(Q)$ are disjoint on the punctured spectrum $\Spec((R/I)_T)-\{T(R/I)_T\}$. However, the above theorem shows that this is not possible. 
\end{proof}
Now let us define the set-up which we need to state and prove our next result.

Let $m$ be a positive integer such that $m \leq n$. Fix a partition $(h_1,\ldots, h_m)$ of $n$, that is, choose positive integers $h_1, \ldots, h_m$ such that $h_1+\ldots +h_m=n$.  
Let $J_i$ be an upper triangular Jordan form of a nilpotent matrix of size $h_i$. For each $h_i$ there are finitely many choices of $J_i$. Let $J=(J_1, \ldots, J_m)$ and let $I_i$ denote the identity matrix of size $h_i$. 

For any $m$-tuple $\underline{\La}=\La_1, \ldots, \La_m$ of distinct elements of $K$, let $J(\underline{\lambda})=J(\La_1, \ldots, \La_m)$ be a matrix such that for all $1\leq i \leq m$, the blocks $\La_iI_i+J_i$ are on the main diagonal. That is, $J(\underline{\La})$ is the direct sum of matrices $\La_iI_i+J_i$. 

Let $\Lambda=\{(\La_1, \ldots, \La_m)  \in \mathbb{A}^m \, | \, \La_i\neq \La_j \text{ for all } 1\leq i\neq j \leq m \}$. It is an open subset of $\mathbb{A}^m$ and therefore is irreducible and has dimension $m$.  Let $$W_J= \{ A\in M_n(K)\, | \, \text{ $A$ is similar to some $J(\La_1, \ldots, \La_m)$ with }$$  $$\La_1, \ldots, \La_m \in K \text{ distinct} \}.$$ 

Let $c_J$ denote the dimension of the set of matrices that commute with $J(\underline{\La})$, for some $\underline{\La}$. This number is independent of the choice of $\underline{\La}$, since $J(\underline{\La})$ commutes with a matrix $A$ if and only if $A$ is a direct sum of matrices $A_i$ such that each $A_i$ has size $h_i$ and $A_i$ commutes with $J_i$. Moreover, $c_J$ is  the dimension of the set of invertible matrices that commute with $J(\underline{\La})$, for some $\underline{\La}$. 

\begin{lem} The dimension of $W_J$ is $n^2 - c_J+m$.   \end{lem}
\begin{proof}
Define a surjective map of algebraic sets 
$$\theta: GL_n(K)\times \Lambda \rightarrow W_J$$
such that $$(U, \La_1, \ldots, \La_m) \rightarrow U^{-1}J(\La_1, \ldots, \La_m)U. $$

Fix $\underline{\La}=(\La_1, \ldots, \La_m)$. Then $$\theta^{-1}(J(\underline{\La}))=\{(U, \underline{\mu}) \in GL_n(K)\times \Lambda \,|\,U^{-1}J(\underline{\mu})U=J(\underline{\La})\}$$  and  it has the dimension of the set $$ \{U \in GL_n(K) \,|\, J(\underline{\La})U=UJ(\underline{\La}) \},$$ that is, it is the set of all invertible matrices commuting with $J(\underline{\La})$. 

Let $M =U^{-1}J(\underline{\La})U\in W_J$ and let $(V, \underline{\mu}) \in \theta^{-1}(M)$. Then $U^{-1}J(\underline{\La})U=V^{-1}J(\underline{\mu})V$  for some $\underline{\mu}$ and $J(\underline{\mu})=(UV^{-1})^{-1}J(\underline{\La})(UV^{-1})$. Hence, $(V, \underline{\mu}) \in \theta^{-1}(J(\underline{\La}))U$.  Therefore, $ \theta^{-1}(J(\underline{\La}))$ and $ \theta^{-1}(M)$ have the same dimension. Since the dimension of $W_J$ is the dimension of $GL_n(K)\times \Lambda$ minus the dimension of a generic fiber $\theta^{-1}(J(\underline{\La}))$, we have that the dimension of $W_J$ is $n^2-c_J+m$.  

Moreover, the set of pairs of matrices $(A,B) \in M_n(K) \times M_n(K)$ such that $A$ and $B$ commute has dimension $(n^2-c_J+m)+c_J=n^2+m\leq n^2+n$.
\end{proof}
\begin{claim} Let $$W=\{(A,B) \in M_n(K) \times M_n(K) | \, [A,B]=0, \, \Delta(A)=0, \, \mathcal{P}(A)=0, \, \mathcal{P}(B)=0\},$$
then there is an injective map $$\Psi:  \mathbb{V}(P+Q, \Delta(X)) \rightarrow W$$
so that $$(A,B)\rightarrow (A,B).$$\end{claim}
\begin{proof}
Let $(A,B) \in \mathbb{V}(Q, \Delta(X)) - \mathbb{V}(P+Q, \Delta(X))$. Then by  Lemma ~\ref{L2Y}, $\mathcal{P}(A)=\mathcal{P}(B)=0$. Therefore, $(A,B) \in W$.  Since $\mathbb{V}(Q)$ is the closure of  $\mathbb{V}(Q) - \mathbb{V}(P+Q)$ we have that $\mathcal{P}(A)=\mathcal{P}(B)=0$ for all $(A,B) \in \mathbb{V}(Q, \Delta(X))$. Hence, $\mathbb{V}(P+Q, \Delta(X)) \subseteq W$. 
\end{proof}

\begin{claim} The dimension of the set $W=\{(A,B) \in M_n(K) \times M_n(K) | \, [A,B]=0, \, \Delta(A)=0, \, \mathcal{P}(A)=0, \, \mathcal{P}(B)=0\}$ is at most $n^2+n-2$.  \end{claim}
\vspace*{-\baselineskip}
\begin{proof}
Let $V= \{(A,B) \in M_n(K) \times M_n(K) | \, [A,B]=0, \, \Delta(A)=0 \}$ and $V_m= \{(A,B) \in V | \,  \text { $A$ has $m$ distinct eigenvalues }\}$. Then we have that $\dim V_m=n^2+m$ and $V=\bigcup_{m=1}^{n-1}V_m$.  Therefore, $\dim V \leq n^2+n-1$. Notice that since $\Delta(A)=0$, $m \leq n-1$. 

Similarly, let $W_m=\{(A,B) \in W | \,  \text { $A$ has $m$ distinct eigenvalues }\}$. Then $W=\bigcup_{m=1}^{n-1}W_m$. For each value of $m$, $W_m \subseteq V_m$. Therefore, the dimension of $W$ is at most $n^2+n-1$. Moreover, $W$ is a closed subset of $V$ defined by the vanishing of $\mathcal{P}(X)$ and $\mathcal{P}(Y)$. To prove the claim we need to show that $\dim W$ cannot be $n^2+n-1$. We do this by showing that $W$ does not contain any component of  $V$ of dimension $n^2+n-1$.  In other words, we show that there are pairs of matrices $(A,B) \in V$ but not in $W$, i.e., either $\mathcal{P}(A)\neq 0$ or $\mathcal{P}(B)\neq 0$. 

Let $A \in M_n(K)$ be a matrix with distinct eigenvalues $\La=\La_1=\La_2, \La_3, \dots, \La_n$.  Then $A$ is similar to a Jordan matrix in two possible forms.

\textbf{Case 1.} $A$ is similar to $J=\left[ \begin{array}{cccccc}
                                                                    \La        & 0    &0& 0&\ldots    & 0\\
                                                                   0 & \La  & 0 & 0& \ldots & 0\\
                                                                     0  &0&\La_3&0     &\ldots &0\\
                                                                     \ldots && \ldots & & \ldots \\
                                                                    0&  0  &0 & 0   &  \ldots  & \La_n \\
                                                                     \end{array} \right]$\\
                                                                     
Take $B=\text{diag}(a_1, \ldots, a_n)$ be a diagonal matrix with distinct entries on the diagonal. Then $[A,B]=0$ and $\mathcal{P}(B)=\prod_{1\leq i<j \leq n}(a_j-a_i)\neq 0$.  \\

\textbf{Case 2.} $A$ is similar to $J=\left[ \begin{array}{cccccc}
                                                                    \La        & 1   &0& 0&\ldots    & 0\\
                                                                   0 & \La  & 0 & 0& \ldots & 0\\
                                                                     0  &0&\La_3&0     &\ldots &0\\
                                                                     \ldots && \ldots & & \ldots \\
                                                                    0&  0  &0 & 0   &  \ldots  & \La_n \\
                                                                     \end{array} \right]$\\

Write $J=\left[ \begin{array}{c|c}
                                 J_0&0\\ \hline 
                                    0&J_1\\
                              \end{array} \right]$, where $J_0=\left[ \begin{array}{ccc}
                                                                                           \La&1&0\\ 
                                                                                           0&\La&0\\
                                                                                           0&0&\La_3
                                                                                           \end{array} \right]$ and $J_1=\left[ \begin{array}{ccccc}
                                                                                                                                            \La_4&0& \ldots &0\\
                                                                                                                                            0&\La_5&\ldots &0\\
                                                                                                                                            & &  \ldots & &\\
                                                                                                                                            0&0&\ldots&\La_n\\
                                                                                                                                            \end{array} \right]$.\\
                                                                                                                                            
\vspace{3mm}
Take  an $n$ by $n$ block-diagonal matrix $U=\left[ \begin{array}{c|c}
                                                                                U_0&0\\ \hline
                                                                                 0&U_1\\
                                                                            \end{array} \right]$ such that $U_0=\left[ \begin{array}{ccc}
                                                                                                                                               0   & 0   & 1  \\ 
                                                                                                                                               1& 1& 0   \\
                                                                                                                                               0&1&1 \\
                                                                                                                                              \end{array} \right]$  and $U_1 \in M_{n-3}(K)$ is the identity matrix. \\
                                                                                                                      
Then   $U^{-1}=\left[ \begin{array}{c|c}
                                    U_0^{-1}&0\\ \hline 
                                     0&U_1^{-1}\\
                                    \end{array} \right]$     and $U^{-1}JU=\left[ \begin{array}{c|c}
                                                                                                                       U_0^{-1}J_0U_0&0\\ \hline 
                                                                                                                       0&J_1\\
                                                                                                                      \end{array} \right]$.\\
                                                                                                                                                                                                                                      
Our goal is to show that $\mathcal{P}(U^{-1}JU)\neq 0$. First, we prove it for the case of 3 by 3 matrices, i.e., for $U_0^{-1}J_0U_0$. 

Observe that $\mathcal{P}(U_0^{-1}J_0U_0)= \mathcal{P}(U_0^{-1}J_0U_0-\La I)=\mathcal{P}(U_0^{-1}(J_0-\La I)U_0)$. 

Denote $M=U^{-1}(J-\La I)U$ and $M_0=U_0^{-1}(J_0-\La I)U_0$.

We have that 
\vspace{3mm}
$$U_0^{-1}=\left[ \begin{array}{ccc}
                                  1   & 1   & -1  \\ 
                                 -1& 0& 1   \\
                                  1&0&0 \\
                                  \end{array} \right]$$ and
\vspace{3mm}                                   
  $$M_0=U_0^{-1}(J_0-\La I)U_0=\left[ \begin{array}{ccc}
                                                                     1   & 1-(\La_3-\La)   & -(\La_3-\La)  \\ 
                                                                      -1& -1+(\La_3-\La)& (\La_3-\La)   \\
                                                                     1 & 1 &0 \\
                                                                     \end{array} \right].$$ 
Moreover,   $$M_0^2=\left[ \begin{array}{ccc}
                                                                     0 & -(\La_3-\La)^2   & -(\La_3-\La)^2  \\ 
                                                                      0&(\La_3-\La)^2& (\La_3-\La)^2   \\
                                                                     0 & 0 &0 \\
                                                                     \end{array} \right].$$ 
                                                                                                                                
In particular, the diagonal $\text{diag}(M_0^i)=(0,(\La_3-\La)^i, 0)$ for all $i \geq 2$.                                                                        
\vspace{3mm}                                                           
Then $$\mathcal{P}(M_0)=\det \left[ \begin{array}{ccc}
                                                         1& 1   &0  \\ 
                                                        1 & -1+(\La_3-\La)& (\La_3-\La)^2   \\
                                                         1 & 0 &0 \\
                                                         \end{array} \right]=(\La_3-\La)^2\neq 0.$$  \\     
 \vspace{4mm}          
Finally, 
\vspace{3mm}
$$\mathcal{P}(M)=\det \left[\begin{array}{cccccc}
                                             1& 1   &0  & 0 & \ldots & 0 \\ 
                                             1 & -1+\La_3-\La& (\La_3-\La)^2 & (\La_3-\La)^3 & \ldots & (\La_3-\La)^{n-1}   \\
                                             1 & 0 &0& 0& \ldots & 0  \\
                                             1& \La_4-\La& (\La_4-\La)^2&( \La_4-\La)^3 & \ldots & (\La_4-\La)^{n-1} \\
                                            \ldots &&\ldots&&\ldots \\
                                             1& \La_n-\La& (\La_n-\La)^2& (\La_n-\La)^3 & \ldots & (\La_n-\La)^{n-1} \\
                                             \end{array} \right] =$$
\vspace{2mm}
 $$=\prod_{i=3}^n(\La_i-\La)^2\prod_{3 \leq i<j\leq n}(\La_j-\La_i)\neq 0.$$                               
The final expression for the determinant is nonzero. Hence $\dim W \leq n^2+n-2$.

Thus we have that $\Delta(X)$ is not in any of the minimal primes of $P+Q$. 
\end{proof}
Now we prove that $P+Q$ has only one minimal prime.
\begin{theorem} Let $P$ and $Q$ be as in Notation ~\ref{NotComm}. Then $ \mathbb{V}(P+Q)$ is irreducible, i.e., $\Rad(P+Q)$ is prime.  \end{theorem}
\vspace*{-\baselineskip}
\begin{proof}
Let $\mathcal{U}$ be a dense open subset in the algebraic set defined by $\mathcal{P}(X)$ as in Lemma ~\ref{L1Y}. Let $A\in M_n(K)$ be such that $\mathcal{P}(A)=0$. Suppose that $A \in \mathcal{U}$. Then by Lemma ~\ref{L1Y} there exists a matrix $B$ such that $(A,B)$ is in the skew-component of the algebraic set of nearly commuting matrices, that is $(A,B) \in \mathbb{V}(Q)$. Let $K[t]$ be a polynomial ring in one independent variable $t$. Fix any $f \in K[t]$. Then $(A,cB+f(A)) \in \mathbb{V}(Q)$ for all $c\in K-\{0\}$. Since $Q$ defines a closed set, we must have that $(A, f(A))  \in \mathbb{V}(Q)$, i.e., when $c=0$ as well.  Since $U$ is a dense subset in $\mathbb{V}(\mathcal{P}(X))$, $(A, f(A)) \in \mathbb{V}(Q)$ for all $A \in \mathbb{V}(\mathcal{P}(X))$. Recall that $f$ was an arbitrary element of $K[t]$. 

Now assume also that $\Delta(A)\neq 0$. Then every matrix $B$ that commutes with $A$ is a polynomial in $A$ of degree at most $n-1$. Thus 
$$ \mathbb{V}(P)_{\Delta(X)}=\{(A, f(A))\, | \, \Delta(A)\neq 0 \text{ and } f  \text{ is a polynomial of degree at most } n-1  \}. $$  Moreover, since $\mathbb{V}(P+Q) \subset \mathbb{V}(P)$, every element of $\mathbb{V}(P+Q)_{\Delta(X)}$ is of the form $(A, f(A))$, where $\mathcal{P}(A)=0$ and $f$ is a polynomial of degree at most $n-1$. 

Identify polynomials $f \in K[t]$ of degree at most $n-1$ with $\mathbb{A}^n$. Then we can consider a map 
$$\mathbb{V}(\mathcal{P}(X))_{\Delta(X)}\times \mathbb{A}^n \rightarrow \mathbb{V}(P+Q)_{\Delta(X)}$$ 
such that $$(A, f) \rightarrow (A, f(A)).$$

Moreover, this map is a bijective morphism. Therefore, $\mathbb{V}(P+Q)_{\Delta(X)}$ is irreducible. If $\mathbb{V}(P+Q)$ is not irreducible, then its nontrivial irreducible decomposition will give us a nontrivial irreducible decomposition of $\mathbb{V}(P+Q)_{\Delta(X)}$, since ${\Delta(X)}$ is not in any of the minimal primes of $P+Q$. Thus the result. 
\end{proof}
\vspace*{-\baselineskip}
\begin{cor} Let $P$ and $Q$ be as in Notation ~\ref{NotComm}. Then, when $n=3$, $P+Q$ is prime. \qed \end{cor}
\vspace{5mm}
\section{The ideal of nearly commuting matrices is a radical ideal}

In this section we prove that $I$ is a radical ideal in all characteristics. We know that $\Rad(I)=P\bigcap Q$ and $I$ is unmixed as the heights of $P$ and $Q$ are equal to $n^2-n$. To prove the result it is sufficient to show that $I$ becomes prime or radical once we localize at $P$ or $Q$. 
\begin{theorem} \label{Rad} The defining ideal of the algebraic set of nearly commuting matrices is radical.\end{theorem}
\begin{proof}
For simplicity of notation, let $\mathcal{P}$ denote $\mathcal{P}(X)$.

We have that $K[X]\bigcap P=(0)$, since every $f \in K[X]\bigcap P$ must vanish when we set $X=Y$. Therefore, $W=K[X]-\{0\}$ is disjoint from $P$ and hence from $I$. Localize at $P$. Then we have an injective homomorphism of $K[X,Y]/I$-modules 
$$(K[X,Y]/I)_P \hookrightarrow (K(X)[Y]/I)_P\cong  (L[Y]/I)_P,$$ 
where $L=K(X)$ and now $I$ is an ideal generated by $n^2-n$ linear equations in the entries of $Y$ with coefficients in $L$. We can always choose at least $n$ variables $y_{ij}, \, (i, j) \in \Lambda$, and write the rest of them as $L$-linear combinations of the chosen ones. Thus $(K[X,Y]/I)_P \hookrightarrow L[y_{ij}]_{(i,j) \in \Lambda}$ and $IK[X,Y]_P$ is prime. 

Next observe that  $K[X] \bigcap Q=(\mathcal{P})$. Clearly, $(\mathcal{P}) \subseteq Q$. To prove the other direction, let $f \in K[X] \bigcap Q$ be nonzero. Then by Lemma ~\ref{L2Y}, $f \in (\mathcal{P})$.  In other words, for all $A \in M_n(K)$ such that $A \in \mathbb{V}(Q)$ and such that there exists a matrix $B$ with the property that $[A,B]$ is nonzero diagonal, then $\mathcal{P}(A)=0$.

Therefore, we have an injective homomorphism of $K[X,Y]/I$-modules
 $$(K[X,Y]/I)_Q \hookrightarrow (V[Y]/I)_Q,$$ 
 where $V=K[X]_{(\mathcal{P})}$ is a discrete valuation domain. 
Then generators of $I$ become linear polynomials in the entries of $Y$ with coefficients in $V$. Let $\mathcal{B}$ be the matrix of coefficients of this linear system such that its rows are indexed by $(i,j)$ for $1\leq i\neq j\leq n$ and columns are indexed by $(h,k)$ for all $1\leq h, k \leq n$. Then for $i\neq h$ and $k\neq j$  $\mathcal{B}$ has an entry $x_{ih}$ in the $(i,k), (h,k)$ spot, has an entry $-x_{kj}$ in the $(i,j), (i,k)$ spot,  and $x_{ii}-x_{jj}$ in the $(i,j)(i,j)$ spot and zero everywhere else. Let $y_1, \ldots, y_{n^2}$ denote the entries of $Y$ such that $y_{(i-1)n+j}=y_{ij}$. In $V[Y]$, $I$ is generated by the entries of the matrix  $$\mathcal{B}\left[ \begin{array}{c}
                                 y_1\\
                                 y_2\\
                                 \ldots\\
                                 y_{n^2}\\
                                 \end{array} \right] .$$

By doing elementary row operations over $V$, we can transform $\mathcal{B}$ into a diagonal matrix $\mathcal{C}$. This gives new generators of $I$. To prove that $IV[Y]$ is radical, it is sufficient to show that the diagonal entries in $\mathcal{C}$ have order at most one in $V$. To this end it reduces to show that $\mathcal{C}$ has rank $n^2-n$ and the ideal generated by the minors of $\mathcal{C}$ of size $n^2-n$  cannot be contained in $\mathcal{P}^2V$. But then it is sufficient to prove this for the original matrix $\mathcal{B}$. Hence it suffices to show:  
\vspace*{-\baselineskip}
\begin{claim}  \begin{itemize}
\item[(1)] The submatrix $\mathcal{B}_0$ of $\mathcal{B}$ obtained from the first $n^2-n$ columns has nonzero determinant in $V$.
\item[(2)] The determinant of $\mathcal{B}_0$ is in $(\mathcal{P})-(\mathcal{P}^2)$.
\end{itemize} \end{claim}
\proof
\item[(1)] It is sufficient to prove the first part of the claim over $K(X)=\text{frac}(V)$, i.e., after we invert $\mathcal{P}$. In this case, since $X$ and $Y$ nearly commute,  they must commute, see Lemma ~\ref{L2Y}.  Moreover, $X$ is a generic matrix, hence its discriminant is nonzero and is not divisible by $\mathcal{P}$. Thus $X$ has distinct eigenvalues and $Y$ is a polynomial in $X$ of degree at most $n-1$. Write $\mathcal{B}=[\mathcal{B}_0|\mathcal{B}_1]$, then our equations become 

$$\mathcal{B}_0\left[ \begin{array}{c}
                                 y_1\\
                                 y_2\\
                                 \ldots\\
                                 y_{n^2-n}\\
                                 \end{array} \right] +\mathcal{B}_1\left[ \begin{array}{c}
                                                                                                           y_{n^2-n+1}\\
                                                                                                           y_{n^2-n+2}\\
                                                                                                             \ldots\\
                                                                                                            y_{n^2}\\
                                                                                             \end{array} \right] =0.$$

Notice that $\mathcal{B}_0$ is invertible if and only if for every choice of the values for $[y_{n^2-n+1}, \ldots, y_{n^2}]$ there is a unique solution for the above equation. 

Furthermore, the bottom rows of $X^0, \, X, \ldots, X^{n-1}$ are linearly independent for a generic matrix $X$. This is true because it even holds for the permutation matrix
$$E=\left[ \begin{array}{cccccc}
0& 1&0& \ldots &0 &0\\
0& 0&1& \ldots &0 &0\\
& \ldots&&\ldots&&\\
0& 0&0& \ldots &0 &1\\
1& 0&0& \ldots &0 &0\\
 \end{array} \right] $$
for which the bottom rows of $E^0, E, \ldots, E^{n-1}$ are the standard basis vectors $e_i$ for $1\leq i \leq n$. 

Thus, given any bottom row $\rho$ of $Y$, there exist $\alpha_0, \ldots, \alpha_{n-1} \in K(X)$ such that $\rho$ equals the bottom row of  $\alpha_0+\alpha_1X+\ldots +\alpha_{n-1}X^{n-1}$. That is,  such a $Y$ is uniquely determined by the entries of its bottom row. Therefore, $\mathcal{B}_0$ is invertible in $K(X)$. 

\item[(2)] First, let us show that $\det \mathcal{B}_0 \in (\mathcal{P})$. For any matrix $A$ in an open dense subset  of the closed set defined by $\mathcal{P}$, there exists a matrix $A'$ such that the commutator $[A,A']$ is a nonzero diagonal matrix, see Lemma~\ref{L1Y} . Hence,  for all $c \in K-\{0\}$ and for all $f\in K[X]$ polynomials of degree at most $n-1$, $(A, cA'+ f(A)) \in \mathbb{V}(I)$. Therefore, the space of solutions of $$\mathcal{B}\cdot \left[ \begin{array}{c}
                                                      y_1\\
                                                      y_2\\
                                                    \ldots\\
                                                   y_{n^2}\\
                                                \end{array} \right] =0$$ has dimension $n+1$, but  we showed  in (1) that it must be $n$. Therefore, the minors of $\mathcal{B}$ must vanish whenever $\mathcal{P}$ vanishes. 

Notice that the degree of the polynomial $\mathcal{P}$ is $n(n-1)/2$, while the degree of $\det \mathcal{B}_0$ is $n(n-1)$. Therefore, to prove part (2) it is sufficient to show that $\det \mathcal{B}_0$ is not a $K$-scalar multiple of $\mathcal{P}^2$. 

Now let us put grading on the entries of $X$ and $Y$. Let $\text{deg } x_{ij}=\text{deg } y_{ij}=i-j$.  Then their products $XY$ and $YX$ and sums have  this property as well: $\text{deg } (XY)_{ij}=i-j$ and $\text{deg } (X+Y)_{ij}=i-j$. Therefore, so does the commutator matrix $XY-YX$.  In fact, any polynomial in  $X$ and $Y$ has this property. Notice that the diagonal entries have degree 0, thus $\mathcal{P}$ has degree 0. 
However, this is not the case for the determinant of the matrix $\mathcal{B}_0$. The nonzero entry corresponding to $(i,j), (k,h)$ has degree $i-j+h-k$. Therefore, if a product of the entries is a nonzero term of the determinant of $\mathcal{B}_0$, then its degree is $\sum_{1\leq i \neq j \leq n}\sum_{1\leq h \leq n }\sum_{1\leq k <n} (i-j+h-k)=n^2(n-1)^2/2\neq 0$ for all $n\geq 2$. Hence $\det \mathcal{B}_0$ cannot be a $K$-scalar multiple of  $\mathcal{P}^2$. That is, when we factor out $\mathcal{P}$ from the minors of $\mathcal{B}$, the remaining expression is not divisible by $\mathcal{P}$. 
\qed

Now we are ready to finish our discussion. Let $\mathcal{C}=[\mathcal{C}_0|\mathcal{C}_1]$ be a matrix that is obtained from $B$ by elementary row transformations so that $\mathcal{C}_0=\text{diag}(c_{ii})_{1\leq i\leq n^2-n}$ is a diagonal matrix. We proved that there exists $1\leq k \leq n^2-n$  with the property that $c_{ii}$ is a unit in $V$ for all $i\neq k$ and $c_{kk}\in(\mathcal{P})-(\mathcal{P}^2)$. Denote $c_{kk}=\alpha \mathcal{P}$, where $\alpha$ is a unit in $V$. 

The ideal $I$ is generated by the following equations
$$\left[ \begin{array}{cccccc}
c_{11}& 0&0&\ldots&0\\
0& c_{22}&0& \ldots &0\\
&& \ldots&&\ldots& \\
0& 0&0& \ldots& c_{n^2-n-1, n^2-n-1}&0\\
0& 0&0& \ldots&0&  c_{n^2-n, n^2-n} \
 \end{array} \right] \left[ \begin{array}{c}
                                            y_1\\
                                           y_2\\
                                          \ldots\\
                                        y_{n^2-n}\\
                                        \end{array} \right]+\mathcal{C}_1\left[ \begin{array}{c}
                                                                                                                y_{n^2-n+1}\\
                                                                                                                y_{n^2-n+2}\\
                                                                                                                       \ldots\\
                                                                                                                      y_{n^2}\\
                                                                                                          \end{array} \right].$$
 Then $V[Y]/I \cong V[y_{n^2-n+1},y_{n^2-n+2}, \ldots, y_{n^2}][y_{kk}]/(\alpha\mathcal{P}-\sum_{j=1}^n (\mathcal{C}_1)_{kj}y_{n^2-n+j})$ is reduced. 
 To show this we consider two cases. If all $(\mathcal{C}_1)_{kj} \in \mathcal{P}$, then the last factor ring is isomorphic to $V[y_{n^2-n+1},y_{n^2-n+2}, \ldots, y_{n^2}][z]/(z\mathcal{P})$. If there is $j$ so that  $(\mathcal{C}_1)_{kj}$ is a unit, then the factor ring is isomorphic to $V[y_i]_{n^2-n+1\leq i \neq j\leq n^2}[y_{kk}]$. In either case, it is reduced.  
  Therefore, since we have an injective map $(R/I)_Q\hookrightarrow (V[Y]/I)_Q$,  $IR_Q$ is radical.
\end{proof}
\vspace{5mm}
\section{Conjectures}

In this section we state conjectures that we have made while doing the research. Many of them appeared as a result of computations performed on a computer algebra program Macaulay2, \cite{GS}.

\begin{conjecture} Let $R$ be as in Notation ~\ref{NotComm}. Then $R/P$, $R/Q$ and $R/(P+Q)$ are $F$-regular. \end{conjecture}

\begin{rem} In the case when $n=2$ the conjecture is true. \end{rem}

The following lemma allows us to reduce the above conjecture to the $F$-regularity of $R/(P+Q)$. 
\begin{lem} Let $R$ be a Noetherian local or $\mathbb{N}$-graded ring of prime characteristic $p>0$ and let $I$ be an ideal (homogeneous in the graded case) generated by a regular sequence.  Let $P$ and $Q$ be ideals of $R$ of the same height such that $P$ and $Q$ are linked via $I=P\bigcap Q$. Let  $R/P$  be Cohen-Macaulay. Suppose that $R/(P+Q)$ is $F$-regular (or equivalently, $F$-rational). Then $R/P$ and $R/Q$ are $F$-regular.   \end{lem}
\begin{proof}
By \cite{PS74}, $R/Q$ is Cohen-Macaulay and has the canonical module isomorphic to $(P+Q)/Q$. Similarly, the canonical module of $R/P$ is $(P+Q)/P$. Then $R/(P+Q)$ is Gorenstein, hence it is $F$-rational if and only if it is $F$-regular. 

Recall that a graded ring $R$ is $F$-regular if and only if $R_{\mm}$ is $F$-regular, \cite{LS99}. Then $R/(P+Q)$ is $F$-rational if and only if its localization at the homogeneous maximal ideal is $F$-rational. Then by applying Corollary 2.9 in \cite{E03} we have that $F$-rationality of  $R/(P+Q)$ implies $F$-regularity of  $R/P$ and $R/Q$.  \end{proof}
Thus if we want to prove that the variety of commuting matrices and the skew component are $F$-regular, it is sufficient to prove the statement for their intersection. Of course we need to know whether $R/P$ is Cohen-Macaulay.

\begin{conjecture} $R/I$ is $F$-pure for all $n$. \end{conjecture}

The above conjecture can be solved by proving the following one. 

\begin{conjecture} Let $\mu=\begin{large}\frac{\prod_{i=1, j=1}^nx_{ij}y_{ij}}{\prod_{i=1}^{n-1}x_{ii}y_{i, n-i+1}\cdot x_{n, n-1}\cdot y_{n-1,1}}\end{large}$. \end{conjecture}

Then $\mu^{p-1}$ is a monomial term of $\prod_{ 1\leq i\neq j \leq n } u_{ij}^p$ with coefficient equal to 1 modulo $p$.  

\begin{rem} The above monomial can be obtained taking the product of all the variables and dividing by the  variables according to the following pattern: denote by $\star$ the variable to be divided out. 

$$X=\left| \begin{array}{cccccc}
                     \star&x_{12} &\ldots &x_{1, n-2}& x_{1, n-1}&x_{1n}   \\
                     x_{21}&\star &\ldots &x_{2, n-2}&x_{2, n-1} &x_{2n}   \\
                     &   \ldots& &&\ldots &    \\
                     x_{n-2, 1}&x_{n-2,2} &\ldots &\star& x_{n-2, n-1}&x_{n-1, n}   \\
                     x_{n-1, 1}&x_{n-1,2} &\ldots&x_{n-1, n-2} &\star& x_{n-1, n}   \\
                     x_{n, 1}&x_{n,2} &\ldots&x_{n, n-2} &\star& x_{n, n}   \\
                        \end{array} \right|, $$
                        \vspace{3mm}
$$Y=\left| \begin{array}{cccccc}
                     y_{11}&y_{12} &y_{1,3}&\ldots & y_{1, n-1}&\star   \\
                     y_{21}&y_{22} &y_{2,3}&\ldots &\star &y_{2n}   \\
                     &   \ldots& &\ldots &    \\
                     y_{n-2, 1}&y_{n-2, 2}&\star& \ldots & y_{n-2, n-1}&y_{n-2, n}\\
                     \star&\star &y_{n-1, 3}&\ldots &y_{n-1, n-1}& y_{n-1, n}   \\
                     y_{n, 1}&y_{n,2} &y_{n, 3}&\ldots &y_{n, n-1}& y_{n, n}   \\                        \end{array} \right| $$
\end{rem}
\vspace{3mm}

\begin{conjecture} Let $X$ be a matrix of indeterminates of size $n$ over a field $K$. Let $\mathcal{P}(X)$ be the irreducible polynomial as in Definition ~\ref{poly}.  Then $K[X]/\mathcal{P}(X)$ is $F$-regular. \end{conjecture}

\begin{conjecture}  The following is a regular sequence on $R/I$ and hence a  part of a system of parameters on $R/J$ and $R/Q$.  \end{conjecture}

$$x_{st}-y_{t, \theta(s,t)},\,  x_{1n}, \,  x_{nn}, \,  x_{11}-y_{2n} $$
\\
for all $1\leq s, t, \leq n$ and where $\theta(s,t)=\left\{\begin{array}{cc}
                                                                                    (s+t) \text{mod } n, & \text{ if } s+t \neq n;\\
                                                                                    n , & \text{ if } s+t=n . \\   
                                                                                    \end{array} \right .$
\begin{rem} The conjecture was verified by using Macaulay2 software when $n=3,4$ over  $K=\mathbb{Q}$ and in some small prime characteristics. \end{rem}

In the case when $n=3$, this is equivalent to the following identifications  of variables in matrices $X$ and $Y$\\
$$X=\left| \begin{array}{cccccc}
                     x_{11}&x_{12} &0   \\
                     x_{21} & x_{22} &x_{22} \\
                     x_{31} &x_{32} &0  \\
                        \end{array} \right|, \quad 
Y=\left| \begin{array}{cccccc}
                     x_{31}&x_{11} &x_{21}   \\
                     x_{22}& x_{32} &x_{12} \\
                     0 &x_{22}  &0  \\
                        \end{array} \right|. $$
\vspace{2mm}
\begin{conjecture} Let $Z \subseteq \{ u_{ij} \, | \, 1\leq i\neq j \leq n \}$ be any subset of cardinality at most $n^2-n-1$. Let $\mathcal{I}_Z$ be the ideal of $R$ generated by the elements of $Z$. Then $R/\mathcal{I}_Z$ is $F$-regular. In particular, $\mathcal{I}_Z$ is a prime ideal. \end{conjecture}
\vspace{5mm}
\section{Acknowledgement}

This paper is the part of the thesis written by the author during her doctorate study at the University of Michigan. The author would like to express an enormous gratitude to Mel Hochster for the support and guidance received  while working on the paper. This work was supported in part by NSF grant DMS-1401384 and the Barbour Scholarship at the University of Michigan. \\
\bibliography{mybib}

\emph{Email address:} zhibek.kadyrsizova@nu.edu.kz\\
Department of Mathematics,\\
School of Science and Technology,\\
Nazarbayev University\\
53 Kabanbay Batyr Ave, Astana 010000, Kazakhstan \\
\end{document}